\documentclass[final]{siamltex}

\usepackage{amssymb}
\usepackage{amsmath}
\usepackage{graphicx}
\usepackage{color}
\usepackage[utf8x]{inputenc} 	
\usepackage{lmodern}     
\usepackage[T1]{fontenc} 
\usepackage{hyperref}
\newcommand{\be}{\begin{equation}}
\newcommand{\ee}{\end{equation}}
\newcommand{\bea}{\begin{eqnarray}}
\newcommand{\eea}{\end{eqnarray}}

\newcommand{\LL}{\mathcal{L}}
\newcommand{\R}{\mathbb{R}}
\newcommand{\RR}{\mathbb{R}}

\newcommand{\Q}{\mathcal{Q}}

\DeclareMathOperator{\diverg}{div}

\newtheorem{remark}{Remark}
\begin{document}
%%%%%%%%%%%%%%
\title{Residual equilibrium schemes for time dependent partial differential equations}
%%%%%%%%%%%%%%

%%%%%%%%%%%%%%

\author{Lorenzo Pareschi\thanks{Mathematics and Computer Science Department, University of
Ferrara, Via Machiavelli 35, 44121 Ferrara, Italy ({\tt lorenzo.pareschi@unife.it}).}
\and
Thomas Rey\thanks{
Laboratoire Paul Painlev\'e, CNRS UMR 8524,
Universit\'e de Lille \& Inria RAPSODI Project,
59655 Villeneuve d'Ascq Cedex, France ({\tt thomas.rey@math.univ-lille1.fr})}
}
\maketitle

%%%%%%%%%%%%%%

	\begin{abstract}
	  Many applications involve partial differential equations which admits nontrivial steady state solutions. The design of schemes which are able to describe correctly these equilibrium states may be challenging for numerical methods, in particular for high order ones. In this paper, inspired by micro-macro decomposition methods for kinetic equations, we present a class of schemes which are capable to preserve the steady state solution and achieve high order accuracy for a class of time dependent partial differential equations including nonlinear diffusion equations and kinetic equations. Extension to systems of conservation laws with source terms are also discussed.
	\end{abstract}

	%%%%%%%%%%%%%%
	\maketitle
	%%%%%%%%%%%%%%
 % \tableofcontents
  %%%%%%%%%%%%%%	

%%%%%%%%%%%%%%%
	
	\section{Introduction}
	  \label{secIntro}
Several applications involve time dependent partial differential equations (PDEs) which admit nontrivial stationary solutions. The design of numerical methods which are capable to describe correctly such steady state solutions may be challenging since they involve the balance between heterogeneous terms like convection, diffusion and other space dependent sources. We refer to \cite{Bouchut04, Gosse13} (and the references therein) for recent surveys on numerical schemes for such problems in the case of balance laws.  

Typical examples, include nonlinear convection-diffusion equations 
	  \begin{equation}
	    \label{eq:prototypeParabolic}
	    \left \{ \begin{aligned}
	      & \frac{\partial u}{\partial t}(t,x) \, = \, \diverg \left( A(x,u(t,x)) + \nabla_x N(u(t,x))\right ), &&  x \in \Omega, \, t > 0,\\
	      & u(0,x) = u_0(x), && x \in \Omega,
	    \end{aligned}\right.
	  \end{equation}
	  for $u = u(t,x) \in \RR$, and hyperbolic balance laws 
\begin{equation}
	    \label{eq:prototypeBalance}
	    \left \{ \begin{aligned}
	      & \frac{\partial u}{\partial t}(t,x) =- \diverg F(u(t,x)) +R(u(t,x)), &&  x \in \Omega, \, t > 0,\\
	      & u(0,x) = u_0(x), && x \in \Omega,
	    \end{aligned}\right.
	  \end{equation}	  
	  for $u = u(t,x) \in \RR^n$ and $n \geq 1$.
Further examples are found when dealing with kinetic equations of the type
	  \begin{equation}
	    \label{eq:prototypeKinetic}
	    \left \{ \begin{aligned}
	      & \frac{\partial f}{\partial t}(t,x,v) =- v\cdot\nabla_x f(t,x,v)\, + \, Q(f)(t,x,v), && x \in \Omega,\, v \in \RR^d, \, t > 0,\\
	      & f(0,x,v) = f_0(x,v) \geq 0, && x \in \Omega,\, v \in \RR^d,
	    \end{aligned}\right.
	  \end{equation}	  	  
	  with $f=f(x,v,t) \geq 0$.
	 In many situations, such systems of equations stabilize onto a large-time behavior which is characterized by an accurate balancing between the different terms appearing on the right hand side of the PDEs. In the sequel, we assume that for the prototype problems \eqref{eq:prototypeParabolic}--\eqref{eq:prototypeKinetic} there is an equilibrium solution $u^{eq}(x)$ or $f^{eq}(x,v)$ such that the right hand side of the corresponding PDE vanish. 
	 	 
	 The construction of numerical schemes which preserves such stationary solutions depends strongly on the particular problem studied. Steady states of convection-diffusion problem are closely related to self-similar solutions of the corresponding diffusion equations \cite{ChangCooper:1970,BuetDellacherie:2010}. Different numerical approaches have been proposed by several authors recently \cite{BessemoulinChatardFilbet:2012,CarrilloChertockHuang:2014,CancesGuichard:2015}. For {balance laws}, there is a very large literature on the so called {well-balanced schemes}, namely schemes which capture the balance between the transport terms and the other ones (we refer to the recent monographies \cite{Bouchut04, Gosse13}). This, in general, is related to the construction of suitable numerical fluxes in space. Classical approaches for kinetic equations include the micro-macro decomposition method \cite{LiuYu:2004, JinShi:2009, LemMieu:2008} and the $\delta f$ method in plasma physics \cite{Dimits93, Denton95}.
	 Finally, for homogeneous {kinetic equations}, one often is interested in schemes which preserve the steady state solution of the collision operator (Gaussian profiles, power laws, etc...). This is related to the velocity discretization and typically is tackled with the aim of {conservative and entropic approximations} \cite{DPR04,DimarcoPareschi15,FilbetPareschiRey:2015}.	  

In this paper we generalize the idea recently introduced in \cite{FilbetPareschiRey:2015} for space homogenous Boltzmann equations to a wide class of PDEs which admit steady states. The main advantage of the approach here presented is its structural simplicity and wide applicability combined with the possibility to achieve very high order of accuracy, which usually is not the case for most approaches present in the literature.

The rest of the paper is organized as follows. In Section \ref{secMicroMacro} we introduce the idea in \cite{FilbetPareschiRey:2015} based on micro-macro decomposition. We show how this can be used to construct steady state preserving schemes for space independent kinetic equations. Next, in Section \ref{secGeneral} we present the residual distribution approach for general PDEs and show its equivalence with micro-macro decomposition for linear problems. Some practical examples of application of the schemes to different prototypes of PDEs are then described in details. In particular, the case of hyperbolic systems with sources is discussed and an analysis of the TVD property of the method is reported in a separate appendix.
In Section \ref{secNumerics} numerical results based on these examples are presented which show the performance of the schemes in different situations. Some conclusions and future developments are reported at the end of the paper.

\section{Micro-macro decompositions}
\label{secMicroMacro}
In this section we recall the approach presented in \cite{FilbetPareschiRey:2015} for kinetic equations of the form (\ref{eq:prototypeKinetic}) in the space homogeneous case
\be
\left \{ \begin{aligned}
	      & \frac{\partial f}{\partial t}(t,v) = Q(f)(t,v), && v \in \RR^d, \, t > 0,\\
	      & f(0,v) = f_0(v), && v \in \RR^d.
	    \end{aligned}\right.,\
	    \label{eq:pkh}
\ee
The approach is based on the classical micro-macro decomposition \cite{LiuYu:2004} and is closely related with other methods presented in the literature, like the works \cite{JinShi:2009, LemMieu:2008} about numerical methods which preserves asymptotic behaviors of some kinetic models or like \cite{CrouseillesLemou:2011, Dimits93, Denton95} for plasma physics.
Reviews of such methods can also be found in \cite{jin:2010portoercole, DimarcoPareschi15}.

\subsection{Linear Fokker-Planck equations}
\label{subLinear}	
As a prototype example we consider the one-dimensional Fokker-Planck equation, where $Q(f)=L(f)$ and $L(f)$ is given by
\be
L(f)=\frac{\partial}{\partial v}\left[(v-u)f+{T}\frac{\partial f}{\partial v}\right],\quad v\in \R.
\label{eq:FP}
\ee
The equation preserves mass, momentum and temperature
\be
\rho=\int_{\R} f\,dv,\quad \rho u = \int_{\R}f\,v\,dv,\quad T=\frac1{\rho}\int_{\R}f\,|v-u|^2\,dv,
\ee
since
\be
\int_{\R} L(f) 
\left(
\begin{array}{c}
  1\\
  v\\
  |v|^2   
\end{array}
\right)
\,dv =0,
\ee
and the equilibrium solutions $f^{eq}$ such that $L(f^{eq})=0$ are characterized by Gaussian distributions 
\be f^{eq}=M(\rho,u,T)=\frac{\rho}{\sqrt{2\pi
T}}\exp\left(\frac{-|v-u|^{2}}{2T}\right). \label{eq:M}\ee
The construction of schemes which preserves the equilibrium states (\ref{eq:M}) is non trivial, the most famous example is given by the {Chang-Cooper method} \cite{ChangCooper:1970}. Extensions of the Chang-Cooper idea to high order and to more general nonlinear Fokker-Planck equations are very difficult \cite{BuetDellacherie:2010}.

The problem is strongly simplified if one introduces a change  of variables based on the \emph{micro-macro decomposition}
\be
g=f-M,
\label{eq:micmac}
\ee  
with $M$ the Maxwellian equilibrium and $g$ such that $\int_{\RR^3} g\,\phi\,dv=0$, $\phi=1,v,|v|^2$. 

Since $L(M)=0$ we have \be L(f)=L(g)+L(M)=L(g),\ee
and therefore, being $M$ time independent, we obtain the equivalent formulation of (\ref{eq:pkh})
\be
\left\{
			\begin{aligned}
		  	\frac{\partial g}{\partial t} &= L(g),\\
		  	f&=M+g.
			\end{aligned}
			\label{eq:micmac_linear}
			\right.
\ee
Now the equilibrium state is $g\equiv 0$ and the construction of numerical approximations of $L$, say $L_h$ for a discretization parameter $h$, such that $L_h(0)=0$ is straightforward even with high order of accuracy. For example, any consistent finite difference approximation to (\ref{eq:FP}) admits $g\equiv 0$ as equilibrium state.  
\begin{remark}
The same linear transformation applies also to the BGK model of the Boltzmann equation, simply considering $Q(f)=Q_{BGK}(f)$ with
\be
Q_{BGK}(f)=\mu(M-f),
\ee
where $\mu$ typically depends on the mass density $\rho$.
\end{remark}
    
    \subsection{The quadratic case}
      \label{subQuad}
The idea just presented applies also to quadratic operators of Boltzmann type, where $Q(f)=\Q(f,f)$ is a symmetric bilinear form 
acting only on $v$ and such that mass, momentum and energy are preserved
\be
\int_{\R} \Q(f,f) 
\left(
\begin{array}{c}
  1\\
  v\\
  |v|^2   
\end{array}
\right)
\,dv =0.
\ee
Here we assume that $\Q$ admits the Maxwellian $M$ in (\ref{eq:M}) a equilibrium state, so that $\Q(M,M)=0$.

Let us consider the decomposition (\ref{eq:micmac}). When
			inserted into the operator $\Q$ gives
			\be
			  \Q(f,f)=\LL(M,g)+\Q(g,g),
			  \label{eq:decom}
			\ee
			where $\LL(M,g)=\Q(g,M)+\Q(M,g)$ is a linear operator.	
Again, since $M$ is independent of $t$, equation \eqref{eq:pkh} now reads
			\be
			\left\{
			\begin{aligned}
		  	\frac{\partial g}{\partial t} &= \LL(M,g)+\Q(g,g),\\
		  	f&=M+g.
			\end{aligned}
			\label{eq:micmac_quadratic}
			\right.
			\ee			
			It is immediate to see that $g\equiv 0$ is an equilibrium state. Therefore the construction of steady state preserving algorithms is based on the construction of approximations of the operator $Q$, say $Q_h$ for a discretization parameter $h$, such that $Q_h(0)=0$. As shown in \cite{FilbetPareschiRey:2015} for the Boltzmann equation this can be achieved rather easily, even with spectral accuracy.  
    
    \section{Residual equilibrium schemes}
      \label{secGeneral}
      In this section we introduce a general approach to construct steady state preserving schemes for time dependent PDEs. We will show at the end of the section how this can be seen as a generalization of the micro-macro decomposition technique described in Section \ref{secMicroMacro}.
      
Suppose we have a differential problem of the form
\be
\frac{\partial u}{\partial t}=G(u),
\label{eq:prot}
\ee
where $G(u)=0$ implies $u(t)=u^{eq}$, with $u^{eq}$ a given {equilibrium state} which may depend on the initial or boundary data of the problem. For example, $G(u)$ can be represented by the right hand side of the prototype problems \eqref{eq:prototypeParabolic}-\eqref{eq:prototypeKinetic}.

Let $G_h$ be an order $q$ approximation of $G(u)$, hereafter called the \emph{underlying method}, which originates the approximated problem
\[
\frac{\partial u_h}{\partial t}=G_h(u_h).
\]
Given the discrete equilibrium state $u^{eq}_h$ we define the \emph{residual equilibrium} as 
\be
{r_h=G_h(u^{eq}_h)},
\label{eq:re}
\ee
note that $r_h=O(h^q)$ and define the new order $q$ approximation 
\[
{{\cal G}_h(u):=G_h(u)-r_h.}
\]
The \emph{residual equilibrium approximation} is then given by
\be
\frac{\partial u_h}{\partial t}={\cal G}_h(u_h).
\label{eq:rea}
\ee
Note that by construction $u_h^{eq}$ is an equilibrium state of (\ref{eq:rea}) since by virtue of (\ref{eq:re}) we have
\[
{\cal G}_h(u_h^{eq})=G_h(u_h^{eq})-r_h=0.
\]

\begin{remark}
The method described above admits several generalizations. For example, it can be extended to capture exact time-dependent solution $u^{s}(t)$ to (\ref{eq:prot}) or self-similar solution such that $u(t)=u^{s}(t)$ for $t\gg 1$. In this case the residual is defined as
\be
r_h=G_h(u^{s}_h(t))-\frac{\partial u^{s}_h(t)}{\partial t},
\label{eq:ss}
\ee
where $u^{s}_h(t)$ is the discrete equivalent of the analytic solution. As a consequence the residual equilibrium approximation (\ref{eq:rea}) preserves exactly the discrete solution $u^{s}_h(t)$. 
%Of course, asymptotically (\ref{eq:ss}) coincides with (\ref{eq:re}) since $u^{s}_h(t)\to u^{eq}$. 
%Similarly it can be extended to quasi steady-state solutions $u^{qs}(t)\to u^{eq}$ as $t\to\infty$. This may be particularly useful when the steady state is not known analytically but it is possible to construct a  
\end{remark}

\subsection{Relation with micro-macro decompositions}    
The above approach can be considered as an extension to arbitrary nonlinear operators of the micro-macro decomposition technique described in Section \ref{secMicroMacro}. 

In fact, let us now consider a semi-discrete scheme applied to the  Fokker-Planck problem (\ref{eq:micmac_linear}) in the micro-macro formulation
\be
\left\{
			\begin{aligned}
		  	\frac{\partial g_h}{\partial t} &= L_h(g_h),\\
		  	f_h&=M_h+g_h,
			\end{aligned}
			\label{eq:micmac_linear2}
			\right.
\ee
where $M_h$ is the discrete equilibrium state. 
We have 
${L_h}(f_h)= {{L_h}(g_h)}+{{L_h}(M_h)}$.
Therefore, since $M_h$ is time independent, in original variables formulation (\ref{eq:micmac_linear2}) is equivalent to 
\be
\frac{\partial f_h}{\partial t}={L_h}(f_h)-L_h(M_h),
\label{eq:reaFP}
\ee   
which is exactly the corresponding residual equilibrium approximation to (\ref{eq:pkh})-(\ref{eq:FP}). It is easy to verify that the same holds true also for the micro-macro formulation (\ref{eq:micmac_quadratic}) in the quadratic case, which in the semi-discrete case is
\be
			\left\{
			\begin{aligned}
		  	\frac{\partial g_h}{\partial t} &= \LL_h(M_h,g_h)+\Q_h(g_h,g_h),\\
		  	f_h&=M_h+g_h.
			\end{aligned}
			\label{eq:micmac_quadratic2}
			\right.
			\ee		
Now we have $\Q_h(f_h,f_h)=\LL_h(M_h,g_h)+\Q_h(g_h,g_h)+\Q_h(M_h,M_h)$ and hence problem (\ref{eq:micmac_quadratic2}) yields the residual distribution formulation
\be
\frac{\partial f_h}{\partial t}={\Q_h}(f_h,f_h)-\Q_h(M_h,M_h).
\label{eq:reaB}
\ee   
More in general, it is easy to verify that the equivalence holds true  for any operator described by a $n$-linear form. Let us finally emphasize that the residual equilibrium formulation is not only more general but also simpler then the micro-macro approach since it avoids the direct discretization of the new terms originated by the micro-macro decomposition (for example the linear operator ${\cal L}_h(M_h,g_h)$ in (\ref{eq:micmac_quadratic2}).      
    
%  \subsection{Convergence Analysis for Linear Problems}
%    \label{subLinear}

%    Discrete dynamical system: convergence towards 0 iff spectrum in %the plane $\{\Re z < 0\}$.

%  \subsection{Nonlinear Problems}
%    \label{subNonlinear}

%    Polynomial nonlinearity only: enough for the numerical examples. 
    
%    Word about nonlocal, convolution based operators?

%  \subsection{Stability Analysis}
%    \label{subStability}
%    For linear problem it is obvious, since thanks to the equivalence with the micro-macro formulation it corresponds to %stability of the underlying method applied to the original problem. For nonlinear problems may be we can consider a %perturbation argument?

\subsection{Local residual equilibrium schemes}
The approach just described can be improved in presence of space-dependent (eventually non smooth) steady states. This is the case for example of hyperbolic balance laws and space non homogeneous kinetic equations. 
 
To illustrate the method, let us consider a one-dimensional balance  law
\be
\frac{\partial u}{\partial t}=-\frac{\partial}{\partial x}F(u)+R(u).
\label{eq:hypcon}
\ee 
%\subsection{Residual equilibrium schemes with limiters}
%\begin{frame}{Local residual equilibrium schemes}
A conservative semi-discrete scheme has the general form 
\be
\frac{\partial u_i}{\partial t} =- \frac{F_{i+1/2}-F_{i-1/2}}{\Delta x}+R_i,
\ee
where $F_{i\pm 1/2}$ are the edge fluxes for the $i$-th cell and 
\[
u_i=\frac1{\Delta x} \int_{x_{i-1/2}}^{x_{i+1/2}} u(x,t)\,dx,\quad R_i=\frac1{\Delta x} \int_{x_{i-1/2}}^{x_{i+1/2}} R(u(x,t))\,dx,
\]
are the cell averages.

Let us denote with $F^{um}_{i\pm 1/2}$ the {numerical fluxes} of the underlying method. 
If we denote by $u^{eq}$ the equilibrium state 
such that
\be
\frac{\partial}{\partial x}F(u^{eq})-R(u^{eq})=0,
\ee
in general, unless the fluxes have been suitably constructed, we have 
\[ 
\frac{F^{eq}_{i+1/2}-F^{eq}_{i-1/2}}{\Delta x} - R^{eq}_i \neq 0,
\]
where $F^{eq}_{i\pm 1/2}$ and $R^{eq}_i$ are defined as $F^{um}_{i\pm 1/2}$ and $R_i$ by replacing $u$ with $u^{eq}$.

Now, we can define the \emph{equilibrium preserving fluxes} ${\cal F}_{i\pm 1/2}$ as
\[
{\cal F}_{i\pm 1/2} = F^{um}_{i\pm 1/2}-F^{eq}_{i\pm 1/2},
\]
%and the equilibrium preserving source \[{\cal R}_i=R_i-R_i^{eq},\]
in order to construct the residual equilibrium semi-discrete scheme 
\be
\frac{\partial u_i}{\partial t} =- \frac{{\cal F}_{i+1/2}-{\cal F}_{i-1/2}}{\Delta x}+(R_i-R^{eq}_i).
\ee
The above scheme is clearly equilibrium preserving, however the modified fluxes may loose some good properties of the original fluxes. These properties, of course, are recovered asymptotically since the flux preserves the exact steady state solution.

The idea is to switch between the two fluxes using a {flux-limiter} based on an equilibrium indicator. We use the standard flux where the solution is far from equilibrium and the well-balanced flux where the solution is close to equilibrium. The resulting flux can be written as
\bea
\nonumber
F_{i\pm 1/2} &=& F^{um}_{i\pm 1/2}-\phi_{i} (F^{um}_{i\pm 1/2}-{\cal F}_{i\pm 1/2})\\[-.25cm]
\label{eq:fl-re}
\\[-.25cm]
\nonumber
&=&F^{um}_{i\pm 1/2}-\phi_{i} F^{eq}_{i\pm 1/2},
\eea
where $\phi_{i}=\phi(r_{i})\in [0,1]$ is the \emph{equilibrium flux limiter} and $r_i$ is a suitable \emph{equilibrium indicator}. We require 
%\[
%r_{i\pm 1/2}=\frac{u_{i\pm 1/2}-u^{eq}_{i\pm 1/2}}{u^{eq}_{i\pm 1/2}}\quad \hbox{with} \quad u_{i\pm 1/2} =\frac{u_{i\pm 1}+u_i}{2},\,\, u^{eq}_{i\pm 1/2} =\frac{u^{eq}_{i\pm 1}+u^{eq}_i}{2}
%\]
%is the \emph{equilibrium indicator}.
that far from equilibrium the limiter is close to zero and the flux is represented by the standard scheme. Similarly, near equilibrium the limiter is close to 1 and the flux is represented by the well balanced scheme. We refer to the Appendix for more details on the role of the equilibrium flux limiter and the choice of the equilibrium indicator.

  \section{Numerical examples}
    \label{secNumerics}
    
    We shall now present in this Section some numerical results illustrating our approach for each of the prototype equations \eqref{eq:prototypeParabolic}, \eqref{eq:prototypeBalance} and \eqref{eq:prototypeKinetic}. 
    Let us start with the linear, isotropic, convection-diffusion case, also known as the linear Fokker-Planck equation \eqref{eq:FP}.
    
    \subsection{The linear Fokker-Planck equation}
      \label{subNumFP}
      
      For the sake of simplicity, we consider the unidimensional case $d = 1$. The domain is then the line $\Omega = [-5,5]$ and we will solve the linear Fokker-Planck equation \eqref{eq:FP}, namely we take in the prototype equation \eqref{eq:prototypeParabolic}
        \[ 
          A(x, r) = (x-u)  r, \qquad N(r) = T r,
        \]
        for $u, T \in \RR$. This equation then preserves the mass $\rho$, the mean momentum $u$ and the temperature $T$:
        \[
          \rho=\int_{\R} f\,dv,\quad \rho u = \int_{\R}f\,v\,dv,\quad T=\frac1{\rho}\int_{\R}f\,|v-u|^2\,dv,
        \]
        The equilibrium solution is given by the Maxwellian distribution
        \[
          u^{eq}(x) = \frac{\rho}{\sqrt{2\pi T}} e^{-(x-u)^2/2T}, \quad \forall x\in \Omega.
        \]        
        The initial condition will be taken as a sum of two gaussians
        \[
          u^{in}(x) = e^{-5(x+2.5)^2} + e^{-5(x-2.5)^2}, \quad \forall x\in \Omega.
        \]

        We will compare the results given by the residual equilibrium scheme \eqref{eq:rea} applied to the upwind (denoted by \textbf{REU} for the residual one, \textbf{SU} for the standard) and central (denoted by \textbf{REC} for the residual one, \textbf{SC} for the standard) approximations of equation \eqref{eq:FP} 
        %(\textcolor{red}{\textbf{Question.}} ref: \cite{leveque:2002}? Presentation of the schemes?)
        , with the Chang-Cooper approximation (\cite{ChangCooper:1970}, denoted by \textbf{CC}).
        In all the simulations, we will take $N_x = 100$ uniformly spaced points with a time step $\Delta t = 1.5\times10^{-4}$.       
        
\begin{figure}
          \begin{center}
          \includegraphics[scale=1]{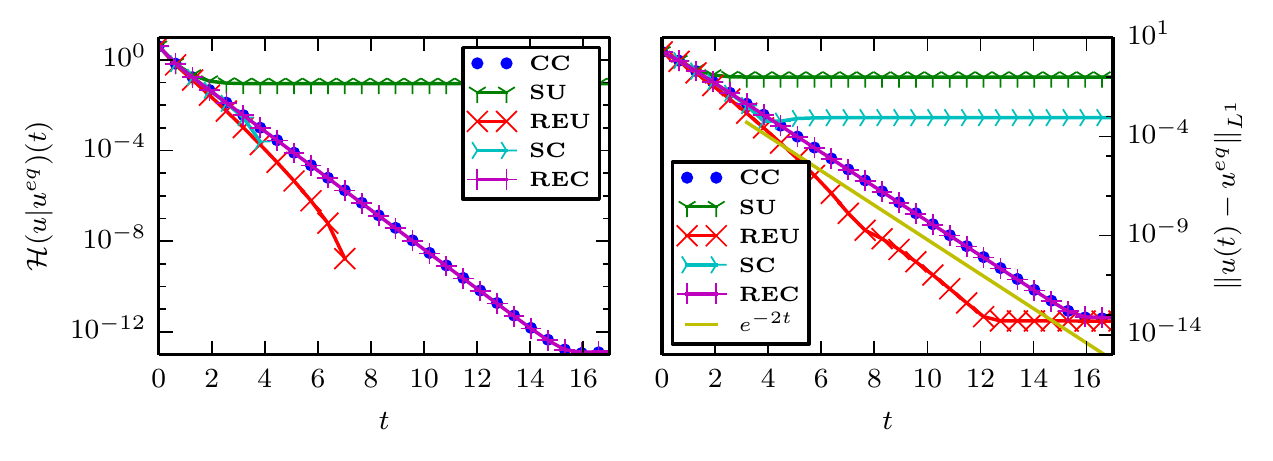}
          \caption{\textbf{Linear Fokker-Planck.} Time evolution of the relative entropy $\mathcal H\left (u|u^{eq}\right )$ (left) and $L^1$ error (right), for the standard upwind (green y's), residual equilibrium upwind (red crosses), standard central (cyan arrows), residual equilibrium central (violet pluses) and Chang-Cooper (blue dots) schemes.}
          \label{fig:LinFP_Error}
          \end{center}
        \end{figure}

        In Figure \ref{fig:LinFP_Error}, we compute the relative entropy
        \[
          \mathcal H\left (u|u^{eq}\right ) := -\int_{\RR} u(x) \log \left (\frac{u(x)}{u^eq(x)}\right )dx,
        \]
        and the $L^1$ error for these five schemes. The quantities are expected to converge exponentially fast toward $0$, with an explicit rate for the later one  \cite{CarrilloToscani:1998}:
        \[ \|u(t) - u^{eq} \|_{L^1} \leq C \, \exp \left ( - 2 t\right ).\]
        As is well known (see \cite{BessemoulinChatardFilbet:2012} for a recent review), the classical upwind (cyan y's) and central schemes saturate (violet arrows), not being able to reach the machine $0$, because of a wrong numerical equilibrium. On the contrary, the Chang-Cooper scheme (blue dots) recovers the correct behavior for both quantities, in particular the correct order of convergence of the $L^1$ norm.
        We observe that the residual equilibrium schemes greatly improve the behavior of the upwind (green crosses) and central schemes (red pluses). 
        First, both quantities decays exponentially toward $0$ in the simulations. More interestingly, the residual equilibrium central schemes gives exactly the same result than the Chang-Cooper scheme, with the great advantage of being very easy to implement and to extend to higher dimension and other models. 
        The correct order of convergence in $L^1$ is recovered only by this latter method, the residual equilibrium upwind not performing  as nicely in this aspect.
        
\begin{figure}
          \includegraphics[scale=1]{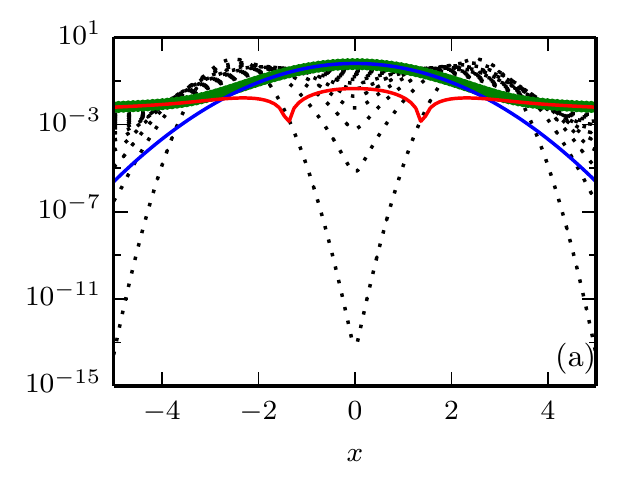} 
          \includegraphics[scale=1]{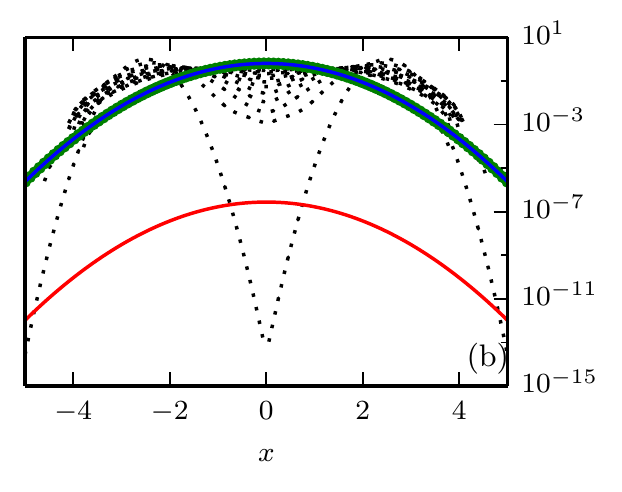} \\
          \includegraphics[scale=1]{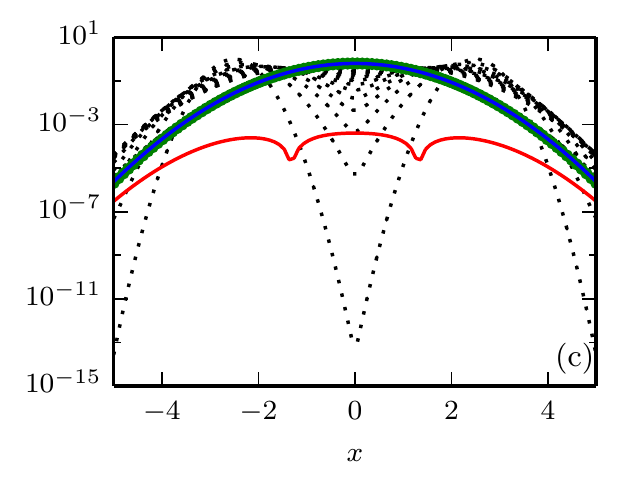} 
          \includegraphics[scale=1]{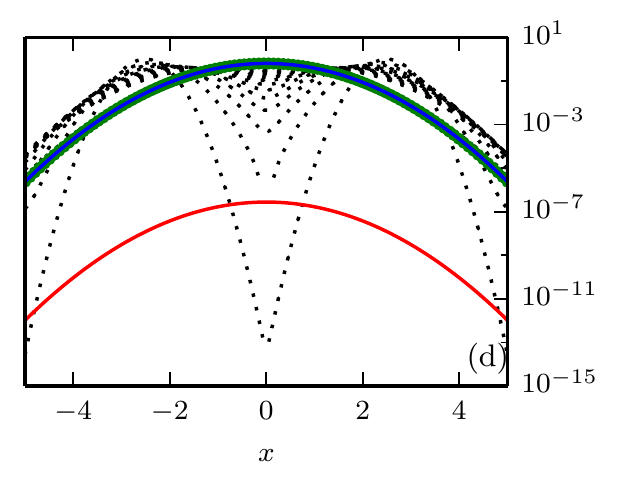}\\
          \includegraphics[scale=1]{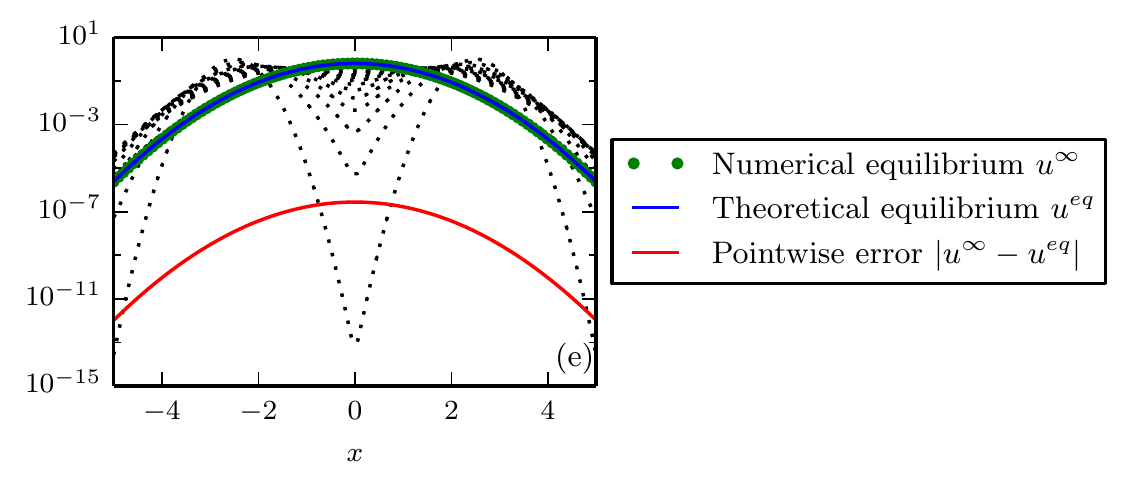}
          \caption{\textbf{Linear Fokker-Planck.} Trend to equilibrium (log scale) of the solutions obtained with the standard upwind (a), residual equilibrium upwind (b), standard central (c), residual equilibrium central (d) and Chang-Cooper (e) schemes.}
          \label{fig:LinFP_Equilibrium}
        \end{figure}

        This behavior can be understood thanks to Figure \ref{fig:LinFP_Equilibrium}, which represents the time evolution of the solution $u(t,x)$ of \eqref{eq:FP} for the $5$ schemes, in semi-logarithmic scale. 
        Although the central scheme (c and d) gives a result very close to the Chang-Cooper (e) one (which is expected according to Figure \ref{fig:LinFP_Equilibrium}), we notice that some oscillations occur for the residual equilibrium upwind scheme in intermediate times (b), as well as some loss of nonnegativity.
        We believe that these oscillations are due to the fact that the residual equilibrium $r_h(t)$ \eqref{eq:re} can be quite large, due to the wrong numerical equilibrium obtained by the classical upwind scheme (a). 
        
%      \paragraph{\textbf{Fermi-Dirac.}}
%    
%        Coefficients:
%        \[ 
%          A(x,r) = xr(1\pm r), \qquad N(r) = r.
%        \]
%        Equilibrium solution:
%        \[
%          u^{eq}(x) = \frac{1}{\beta e^{x^2/2\sigma} \pm 1}.
%        \]
     
    \subsection{The porous medium equation}
      \label{subNumPorMed}
      
      We are now interested in the nonlinear porous medium equation. It is obtained by taking the prototype equation \eqref{eq:prototypeParabolic} with the coefficients
      \[ 
        A(v, r) = v \qquad N(r) = |r|^m.
      \]
      The equation obtained (in self-similar form) can be read as: 
      \begin{equation}
        \label{eq:porousmedia}
        \partial_t f = \nabla_v \cdot \left ( v \, u \right) + \Delta_v\left (u^m\right ).
      \end{equation}
      It is well known (see \cite{CarrilloToscani:2000} for a comprehensive review) that the equilibrium solution of this equation are given by the so called Barrenblatt-Pattle distribution 
      \[
        u^{eq}(v) = \left (C - \frac{m-1}{2m}|v|^2\right )_+^{1/(m-1)}, \quad \forall x \in \Omega,
      \]
      where the constant $C$ depends on the initial condition, and insure mass preservation.
      The solution to \eqref{eq:porousmedia} converges exponentially fasts toward this equilibrium profile.
        
      In all our simulations, we will take $m = 5$ in the equation, and consider  $d=2$, $\Omega = [-10,10] \times [-10,10]$ as the computational domain, with an explicit solver needing a parabolic time step. The initial datum will be chosen as:
      \[
        u^{in}(x) = |x|^2 e^{-|x|^2/2}.
      \]
      We will take $N_x = N_y = 64$ points in each of the space dimension and we will compare the numerical results obtained with a first order non equilibrium-preserving upwind scheme (see \cite{ChainaisPeng:2004} for an analysis of this method, \textbf{SU}) with its residual equilibrium counterpart (\textbf{REU}) and some second order, equilibrium preserving methods: the fully upwind (see \cite{BessemoulinChatardFilbet:2012}, \textbf{FU}) and the Scharfetter-Gummel schemes (see \cite{LazarovMishevVassilevsk:1996}, \textbf{SG}).

        \begin{figure}
          \begin{center}
          \includegraphics[scale=1]{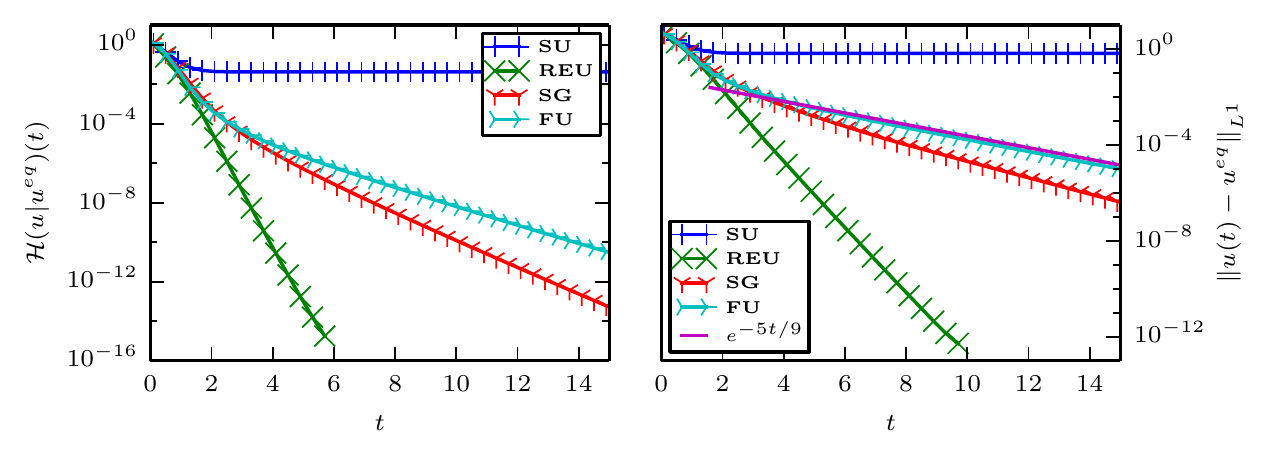}
          \caption{\textbf{Porous medium.} Time evolution of the relative entropy (left) and $L^1$ error (right), for the classical upwind (blue pluses), residual equilibrium upwind (green crosses), fully upwind (cyan arrows) and Scharfetter-Gummel (red y's) schemes.}
          \label{fig:PM_Error}
          \end{center}
        \end{figure}

%      Initial condition. Sum of Gaussians:
%      \[
%        f(v) = \left\{ \begin{aligned}
%                   & e^{\frac{-1}{6-|v-c_0|}} && |v-c_0| < 6, \\
%                   & e^{\frac{-1}{6-|v-c_1|}} && |v-c_1| < 6
%               \end{aligned} \right.
%      \]
%      with $c_0 = (2,-2)$ and $c_1=(-2,2)$.
%      
%        
        In Figure \ref{fig:PM_Error}, we compute the relative entropy
        \[
          \mathcal H\left (u|u^{eq}\right ) := -\int_{\Omega} \left ( \left ( u - u^{eq}\right ) + \frac{2}{m-1} \left (u^m - (u^{eq})^m\right )\right ) dx,
        \]
        and the $L^1$ error for these four schemes. The quantities are expected to converge exponentially fast toward $0$, and we can moreover compute the decay rate for the $L^1$ error, namely according to \cite{CarrilloToscani:2000}
        \[ \|u(t) - u^{eq} \|_{L^1} \leq C \, \exp \left ( - \frac{d(m-1) + 2}{(d+2)m - d} t\right ).\]
        As is expected again, the classical upwind (blue crosses) saturates, not being able to reach the machine $0$, because of a wrong numerical equilibrium. On the contrary, the other three schemes behave properly regarding to their equilibrium preserving properties: they all succeed to capture the Barrenblatt-Pattle distribution.
        Nevertheless, we notice that only the fully upwind scheme (cyan arrows) from \cite{BessemoulinChatardFilbet:2012} is able to capture the correct decay rate. Our new residual equilibrium scheme converges too fast toward the equilibrium. According to the previous section, we believe that this could be improved by switching to another numerical scheme, such as basic central discretization.

       \begin{figure}
         \includegraphics[scale=1]{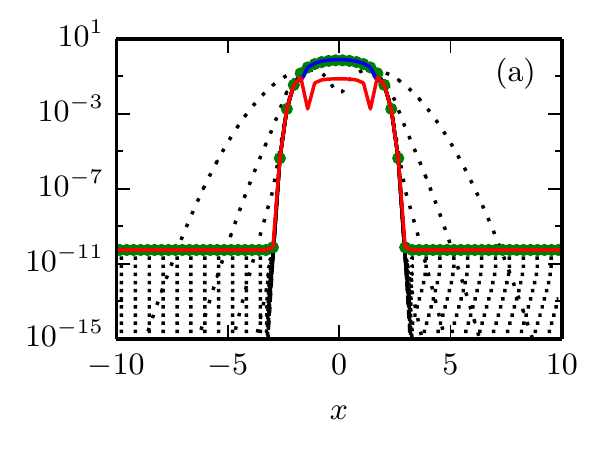} 
         \includegraphics[scale=1]{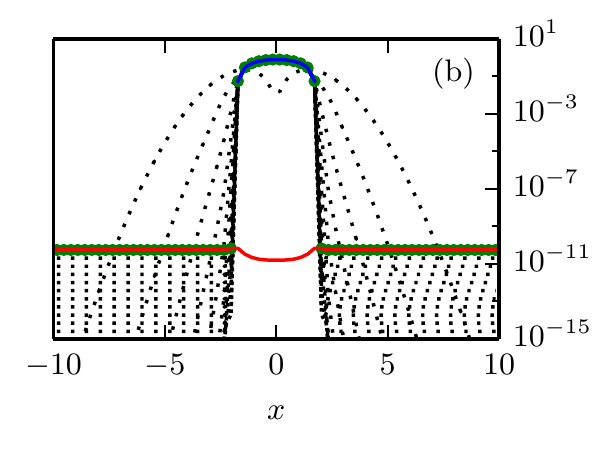}  \\
         \includegraphics[scale=1]{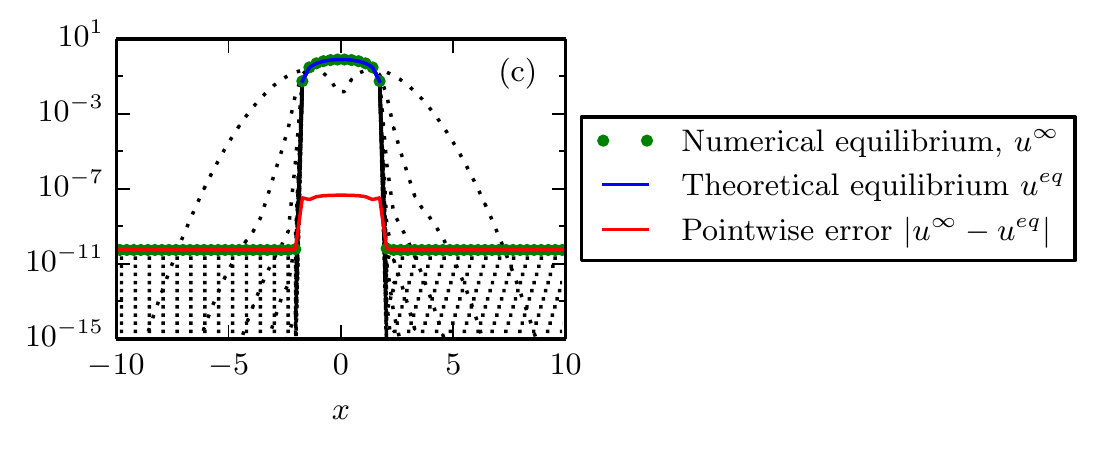} 
         \caption{\textbf{Porous medium.} Trend to equilibrium (log scale) of the solutions obtained with the standard upwind (a), residual equilibrium upwind (b), and fully upwind (c) schemes, at $y = 0$.}
         \label{fig:PM_equilibrium}
      \end{figure}
      
\begin{figure}
        \begin{center}
        \includegraphics[scale=1.1]{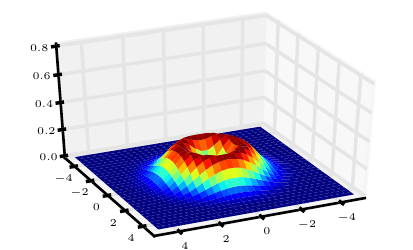}\hskip -.55cm 
        \includegraphics[scale=1.1]{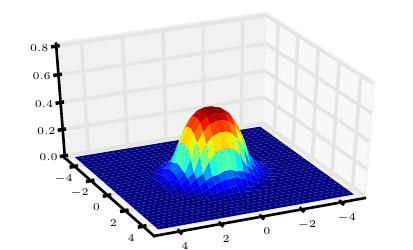}\hskip -.55cm 
        \includegraphics[scale=1.1]{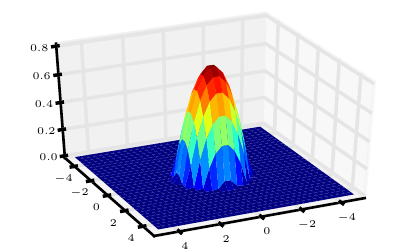} 
        \caption{\textbf{Porous medium.} Solution obtained with the residual equilibrium upwind scheme, for $t=0$, $0.5$ and $20$.}
        \label{fig:PM_evol3d}
        \end{center}
      \end{figure}

      These observations are confirmed by Figure \ref{fig:PM_equilibrium}, which presents in log scale the time evolution of the solution $u(t,x,0)$ for $x \in [-10,10]$ to \eqref{eq:porousmedia} given by the standard upwind (a), residual equilibrium upwind (b) and fully upwind (c, for a reference solution) schemes. 
      We observe that even if the rate of convergence is not correct, the  use of a residual equilibrium approach allows the upwind scheme to compute very nicely the Barenblatt-Pattle distribution (see also Figure \ref{fig:PM_evol3d}), the $L^\infty$ error being of the same order of magnitude than the one obtained with the fully upwind scheme.
      
    \subsection{Boltzmann Equation}
      \label{subBoltzmann}

			We now consider the prototype equation \eqref{eq:prototypeKinetic}: we present a numerical example for the space homogeneous Boltzmann equation in $2$ velocity dimensions, with Maxwell molecules
			\begin{equation}
			  \label{eq:BoltmannHomog}
			  \partial_t f = \int_{\R^{2}} \int_{\mathbb{S}^{1}}   \left[ f'_* f' - f_* f \right] \, d\omega \, dv_*,
			\end{equation}
			where we used the shorthands $f = f(v)$, $f_* = f(v_*)$, $f ^{'} = f(v')$, $f_* ^{'} = f(v_* ^{'})$. 
			The velocities of the colliding pairs $(v,v_*)$ and $(v',v'_*)$ are related by
			\begin{equation*}%\label{eq:rel:vit}
			  v' = \frac{v+v_*}{2} + \frac{|v-v_*|}{2} \sigma, \qquad
			  v'_* = \frac{v+v^*}{2} - \frac{|v-v_*|}{2} \sigma\nonumber.
			\end{equation*}
			As for the linear Fokker-Planck equation \eqref{eq:FP}, the mass $\rho$, the mean momentum $u$ and the temperature $T$ are preserved, and the equilibrium solution is given by the Maxwellian distribution
      \[
        f^{eq}(x) = \frac{\rho}{\sqrt{2\pi T}} e^{-(x-u)^2/2T}, \quad \forall x\in \Omega.
      \]        
 
      We compare the fast spectral method (denote by \textbf{FS}) of \cite{MoPa:2006} with its residual preserving (\textbf{REFS}) counterpart.
			For this, we use an exact solution of the homogeneous Boltzmann equation, the so called Bobylev-Krook-Wu solution \cite{Bobylev:75,KrookWu:1977}. 
			It is given by
			\begin{equation*}
		  	f_{BKW}(t,v) = \frac{\exp(-v^2/2S)}{2\pi S^2} \,\left[2\,S-1+\frac{1-S}{2 \,S}\,v^2 \right]
			\end{equation*}
			with $S = S(t) = 1-\exp(-t/8)/2$. 
			For the velocity discretization, we choose $\Omega = [-8,8]^2$, $N_{v_x} = N_{v_y} = 64$ points and $M = 8$ angular discretizations. Because the problem is not steep, we can take $\Delta t = 0.01$.

      \begin{figure}
        \includegraphics[scale=1.]{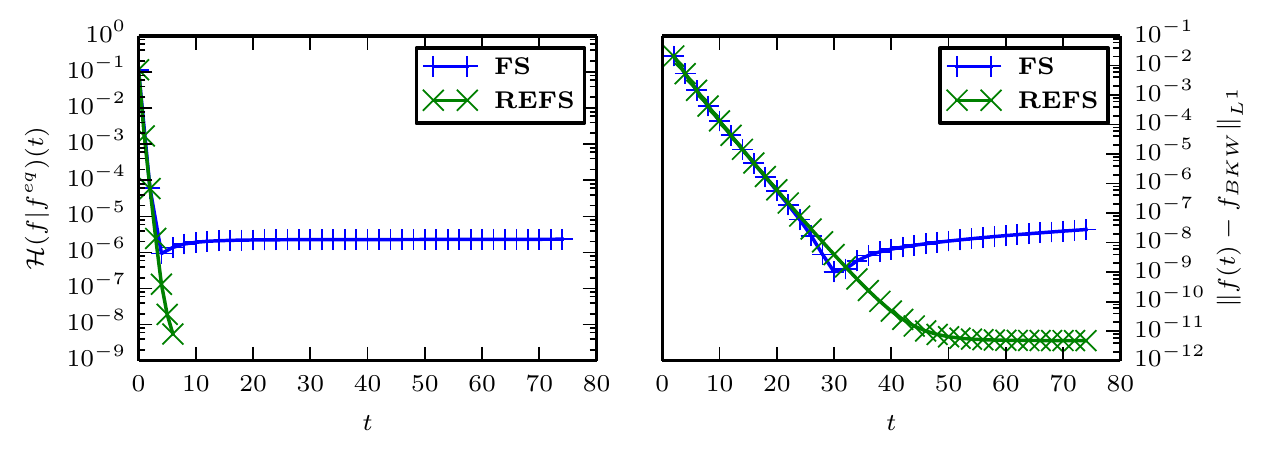}
        \caption{\textbf{Boltzmann equation.} Time evolution of the relative entropy $\mathcal H(f|g)$ (left) and the $L^1$ error (right) with respect to the exact Bobylev-Krook-Wu solution, for the residual equilibrium (green, \textbf{REFS}) and the classical fast spectral method (blue, \textbf{FS}).}
        \label{fig:BoltzmannError}
      \end{figure}

      Figure \ref{fig:BoltzmannError} presents a comparison between both methods for the relative entropy of the solution $f$ with respect to the equilibrium $f^{eq}$
		  \[
		    \mathcal H(f|g)(t) := \int_{\R^d} f(t,v) \log \left (\frac{f(t,v)}{M(v)} \right ),
		  \]
      and the absolute $L^1$ error between the numerical solution $f(t,v)$ and the exact one $f_{BKW}(t,v)$.
      As observed in \cite{FilbetPareschiRey:2015}, for both quantities, the behavior of the residual equilibrium method is better than the classical one. 
      In particular, the new method achieves a nice monotonous decay of the relative entropy, without the large time increase of the classical spectral method. This is due to the fact that the equilibrium of this latter methods are constants \cite{FilbetMouhot:2011}.

      \begin{figure}
        \begin{center}
        \includegraphics[scale=1.1]{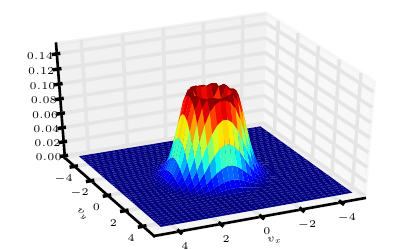}\hskip -.55cm
        \includegraphics[scale=1.1]{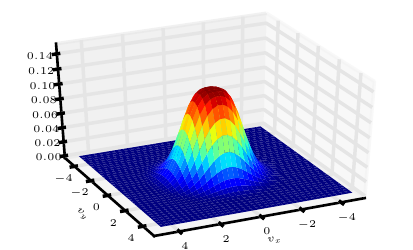}\hskip -.55cm 
        \includegraphics[scale=1.1]{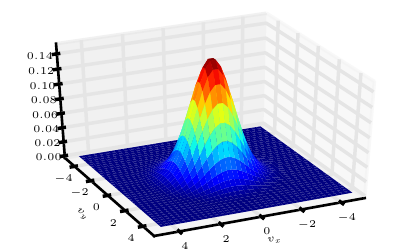} 
        \caption{\textbf{Boltzmann equation.} Solution obtained with the residual equilibrium fast spectral scheme, for $t=0$, $0.5$ and $50$.}
        \label{fig:Boltzmann_evol3d}
        \end{center}
      \end{figure} 
                
      Finally, as a last observation, we noticed (see also Figure \ref{fig:Boltzmann_evol3d}) that even if the positivity of the solution is not preserved by the new method, the number of nonpositive cells seems to be lower at a fixed time for the steady-state preserving method, compared to the classical one.
                        
    \subsection{Shallow Water}
      \label{subShallow}
      
      Let us now consider a simple nonlinear system of balance laws, the so-called one-dimensional shallow water equations with topography
      \begin{equation}
        \left \{ \begin{aligned}
          & \partial_t h + \partial_x(h v) = 0, \\
          & \partial_t(h v) + \partial_x \left (h v^2 + \frac12g h^2\right ) = - g h \partial_x B.
        \end{aligned} \right.
        \label{eq:shallowwater1d}
      \end{equation}
      The unknowns $(h,v)$ denote respectively the water height above the bottom and the mean velocity of the fluid, $g$ is the gravitational constant, and $B(x)$ is the bottom topography (so that $B(x) + h(t,x)$ denotes the top surface). 
      It fits in the framework of the prototype equation \eqref{eq:prototypeBalance} with $n=2$ and
      \[
        u =    \begin{pmatrix} h \\ hv \end{pmatrix} =: \begin{pmatrix} u_1 \\ u_2 \end{pmatrix},  \qquad 
        F(u) = \begin{pmatrix} u_2 \\ u_2^2/u_1 + g u_1^2/2 \end{pmatrix}, \qquad 
        R(u) = \begin{pmatrix} 0 \\ -g u_1 \partial_x B \end{pmatrix}.
      \]
      
      A presentation of some classical approach to solve numerically this system, in both one and two space dimension, can be found in e.g. \cite{Wen:2006}. 
      Due to the strictly hyperbolic nature of the equation \eqref{eq:shallowwater1d} (the eigenvalues $u \pm \sqrt{g h}$ being distinct as soon as $h \neq 0$, namely as soon as there is no dry state), we shall use in all our numerical experiments a second order Lax-Friedrichs scheme with Van Leer's flux limiters \cite{VanLeer:1977a}.
      The numerical parameters will always verify $\Delta t/ \Delta x = 10$.
      We will present two test cases, a simple perturbation of a lake at rest and then a more complex transient, transcritical flow with shock. We chose these tests because it is possible in both cases to compute analytical steady states.
      
      \paragraph{\textbf{Small perturbation of a lake at rest}}
          
        This well known example has been taken from \cite{Leveque:1998}. We assume that there is initially no motion of the fluid and that the surface is flat. 
        This is clearly a steady state of \eqref{eq:shallowwater1d}
        \[ h^{eq}(x), v^{eq}(x)) = (1 - B(x),0) \quad \forall x \in \Omega,\]
        that we will perturb in a small zone:
        \[
          (h^{in}(x), v^{in}(x)) = \left (1+\varepsilon(x) - B(x), 0\right ),
        \]
        with $\varepsilon(x) = 0.1$ if $0.1 < x < 0.2$ and $0$ elsewhere.
        The bottom topography is chosen as the single bump
        \begin{equation}
          \label{eq:topography1}
          B(x) =  \left\{ \begin{aligned}& 0.25\left (\cos(\pi(x - 0.5)/0.1 +1 \right ), && |x - 0.5| < 0.1 \\
                                         & 0 && \text{otherwise}
                 \end{aligned}\right. 
        \end{equation}
        To deal with an infinite lake, we assume no flux at the boundaries, and take $\Omega = [0,1]$. 
        We perform numerical simulations with the classical Lax-Friedrichs (denoted by \textbf{LF}) and the residual equilibrium Lax-Frierichs (\textbf{RELF}) schemes.

        \begin{figure}
          \begin{center}
          \begin{tabular}{c}
            \includegraphics[scale=1]{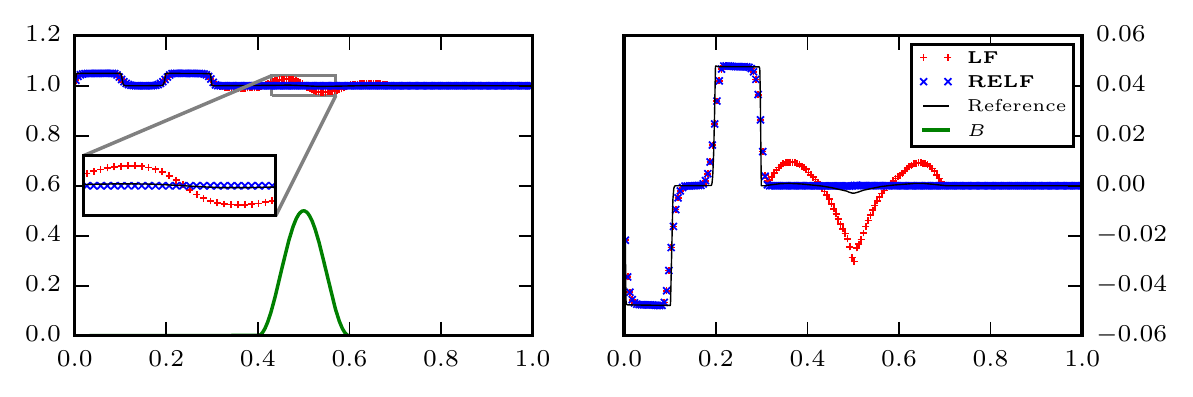} \\
            \includegraphics[scale=1]{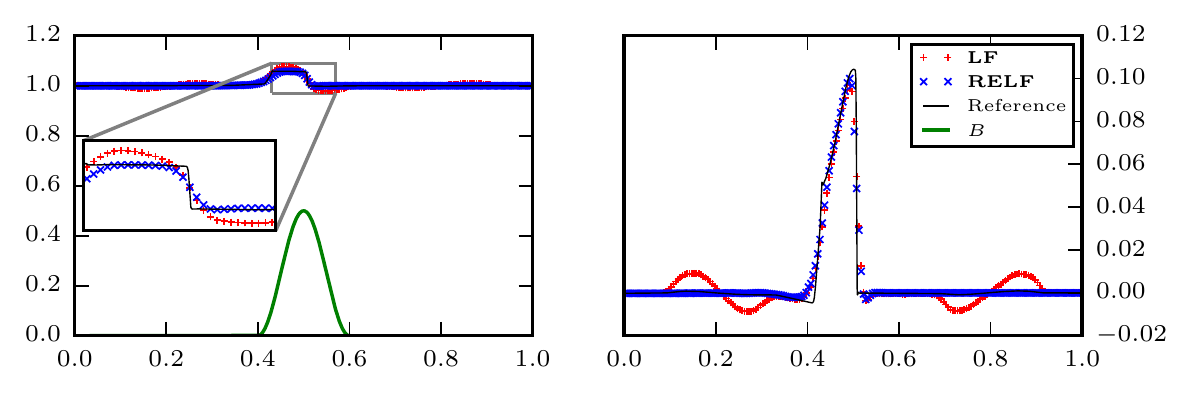} \\
            \includegraphics[scale=1]{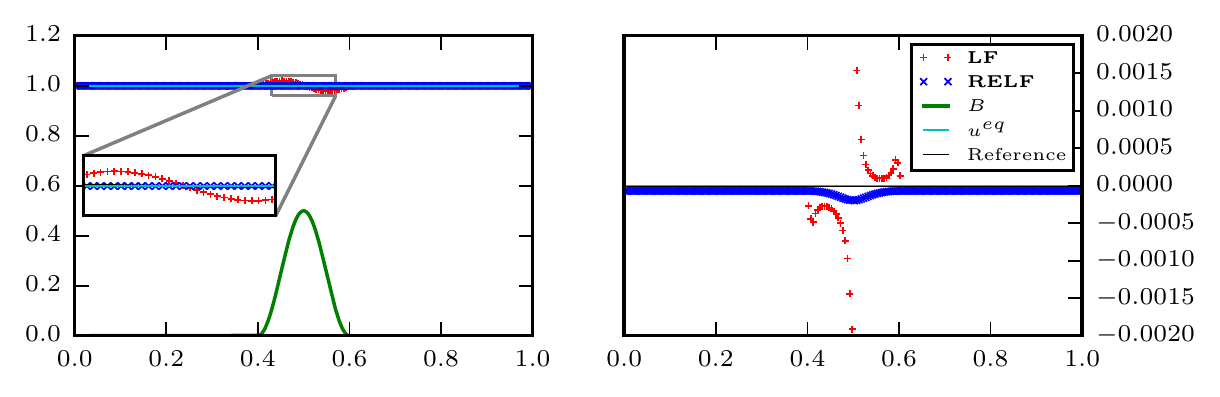}
          \end{tabular}
          \end{center}
          \caption{\textbf{Shallow water: lake at rest.} Evolution of the water height (left) and water velocity (right) at times $t=0.03$, $0.1$ and $1$, for the classical (red dots) and equilibrium preserving (blue dashes) Lax-Friedrichs schemes. The green line is the topography and the cyan one, the equilibrium solution.}
          \label{fig:SW_Lake_Evol}
        \end{figure}
        
        Figure \ref{fig:SW_Lake_Evol} show the time evolution of the solution obtained with both classical (\textbf{LF}) and residual equilibrium (\textbf{RELF}) Lax-Friedrichs schemes, with $N_x = 200$ points in space. The so-called reference solution is given by the \textbf{LF} solver used with $N_x = 1600$ points in space. 
        We notice that the new residual equilibrium solver is much better at capturing both large and (more surprisingly) short time behaviors of the solution. It is clear by  looking at the velocity profiles that it also avoids spurious oscillations induced by the flux limited Lax-Friedrichs scheme.

        \begin{figure}
          \begin{center}
            \includegraphics[scale=.8]{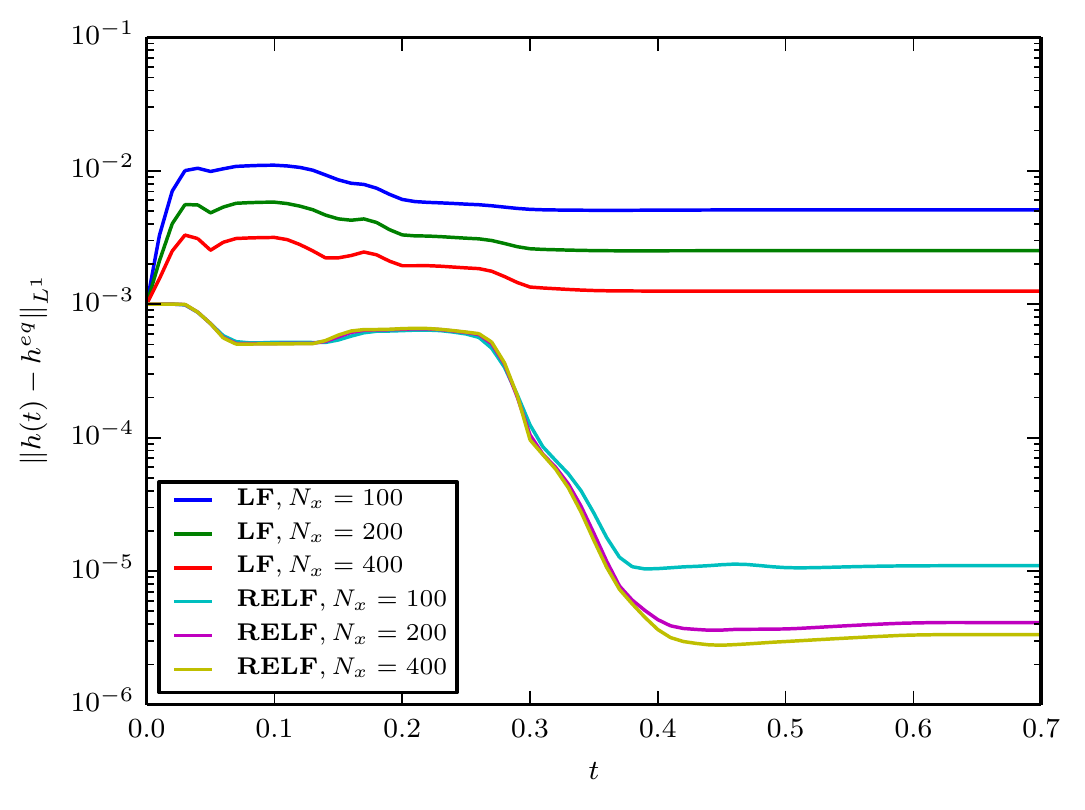}
          \end{center}
          \caption{\textbf{Shallow water: lake at rest.} Time evolution of the $L^1$ error of the surface water $h$, with respect to the equilibrium profile, for the standard (\textbf{LF}) and residual equilibrium (\textbf{RELF}) Lax-Friedrichs schemes, with different mesh sizes.}
          \label{fig:SW_Lake_error}
        \end{figure}
        
        Figure \ref{fig:SW_Lake_error} show the time evolution of the $L^1$ error of the surface water with respect to the equilibrium profile.
        The model is not dissipative, in the sense that there is no global functionals $\mathcal H$ and $\mathcal I \geq 0$ verifying the following type of relations:
        \[
          \frac{\partial \mathcal H}{\partial t} [u](t) \leq -\mathcal I[u](t),
        \]
        for $u$ solution to \eqref{eq:shallowwater1d}. This is usually the case for models \eqref{eq:prototypeParabolic} and \eqref{eq:prototypeKinetic}. Then, one cannot expect monotone (exponential or polynomial) convergence of $u$ toward the equilibrium solution $u^{eq}$, as is well known and was observed numerically in the other models.
        Nevertheless, we can see in Figure \ref{fig:SW_Lake_error} that the solutions obtained with the  \textbf{RELF} scheme converges toward the equilibrium, which is not the case of the \textbf{LF} scheme.

      \paragraph{\textbf{Transcritical flow with shock}}
          
        We now consider the case of a transcritical flow, which will generate a non-smooth equilibrium: a shock will form after the obstacle in the bottom and remain stable in time. 
        This test is taken from \cite{GoutalMaurel:1997}, where the bottom function is not $\mathcal{C}^1$:
        \begin{equation}
          B(x) =  \left\{ \begin{aligned}& \left(0.2 − 0.05(x − 10)^2 \right ), && |x - 10| < 2, \\
                                           & 0 && \text{otherwise}.
                   \end{aligned}\right.
        \end{equation}
        The initial condition is given by the constant profiles
        \[
            (h^{in}(x), v^{in}(x)) = \left (0.33 - B(x), 0\right ), \quad \forall x \in \Omega,
        \]
        where the computational domain $\Omega = [0,25]$. We suppose an inflow $h v = 0.18$ on the left boundary $x=0$, and an outflow $h = 0.33$ on the right boundary $x=25$.
        
        The flow is subcritical (the scaled Froude number $Fr := |v|/\sqrt{g h}$ being less than $1$) before the hump, and supercritical after it: a shock is then generated around the obstacle.
        To compute the equilibrium profile, one has first to locate the shock using the Rankine-Hugoniot's relations, and then use the Bernoulli relation (link between the topography and the water height) outside the shock. 
        The equilibrium profile is then given as a solution of a nonlinear equation, using a Newton solver.

        \begin{figure}
          \begin{center}
          \begin{tabular}{c}
            \includegraphics[scale=1]{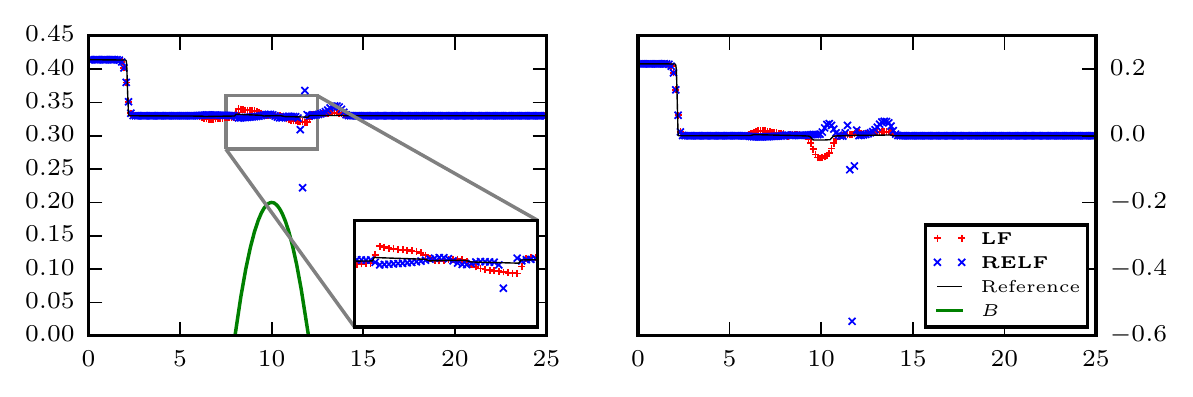} \\
            \includegraphics[scale=1]{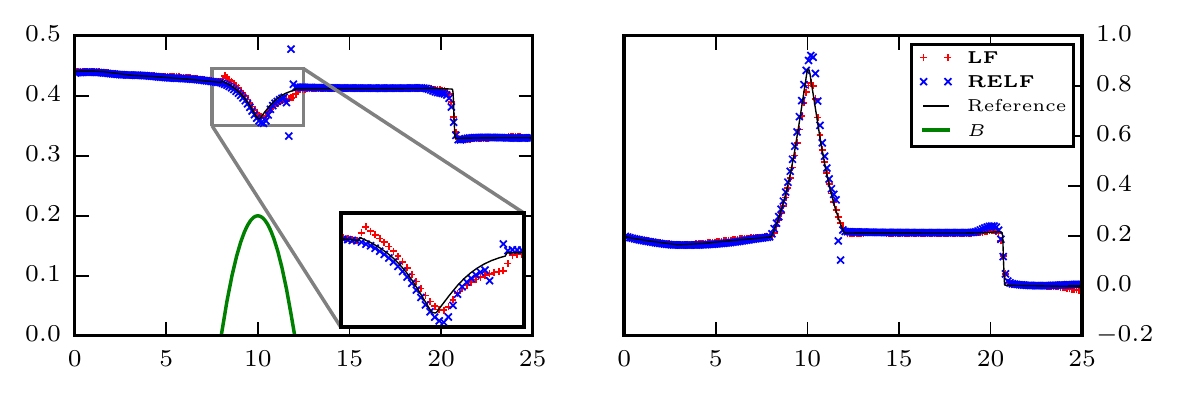} \\
            \includegraphics[scale=1]{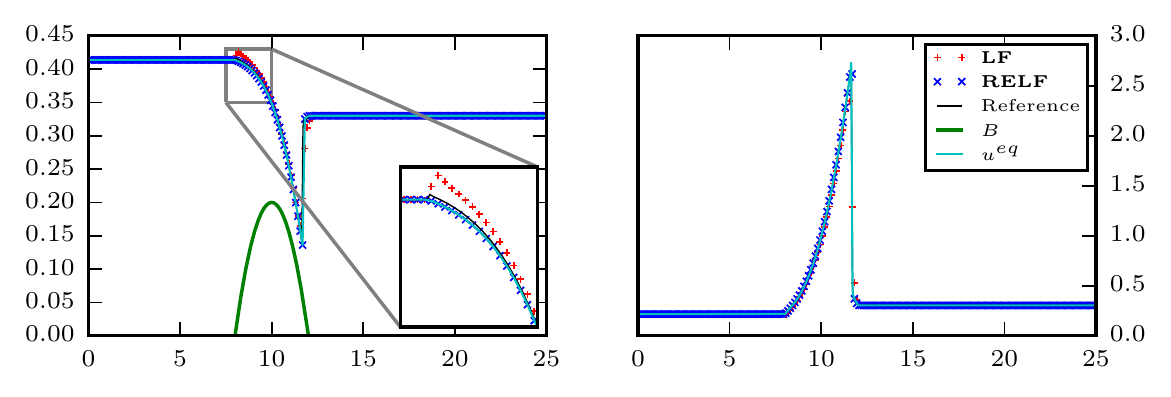}
          \end{tabular}
          \end{center}
          \caption{\textbf{Shallow water: transcritical flow.} Evolution of the water height (left) and Froude number $Fr= |v|/\sqrt{g h}$ (right) at times $t=1$, $10$ and $500$, for the classical (red dots) and equilibrium preserving (blue dashes) Lax-Friedrichs schemes. The green line represents the topography, the black dashes the reference solution and the cyan line, the equilibrium.}
          \label{fig:SW_transcritical_Evol}
        \end{figure}    
        
        Figure \ref{fig:SW_transcritical_Evol} presents the evolution of the solutions obtained with both \textbf{LF} and \textbf{RELF} schemes, at times $t=1$, $10$ and $500$ (equilibrium), for $N_x = 200$ points in space. The reference solution has been obtained by the \textbf{LF} solver using $N_x = 1600$ points.
         Although the new \textbf{RELF} method avoids some of the oscillations obtained with the \textbf{LF} one, we notice that the solution obtained with this new scheme generates a spurious $\mathcal O(1)$-oscillation at the position of the equilibrium shock. 
         Some more experiments show that this is independent of the grid size, and seems actually due to the fact that the equilibrium solution being not smooth, the residual equilibrium \eqref{eq:re} is big and perturbs locally the behavior of the scheme.

        \begin{figure}
          \begin{center}
          \begin{tabular}{c}
            \includegraphics[scale=1]{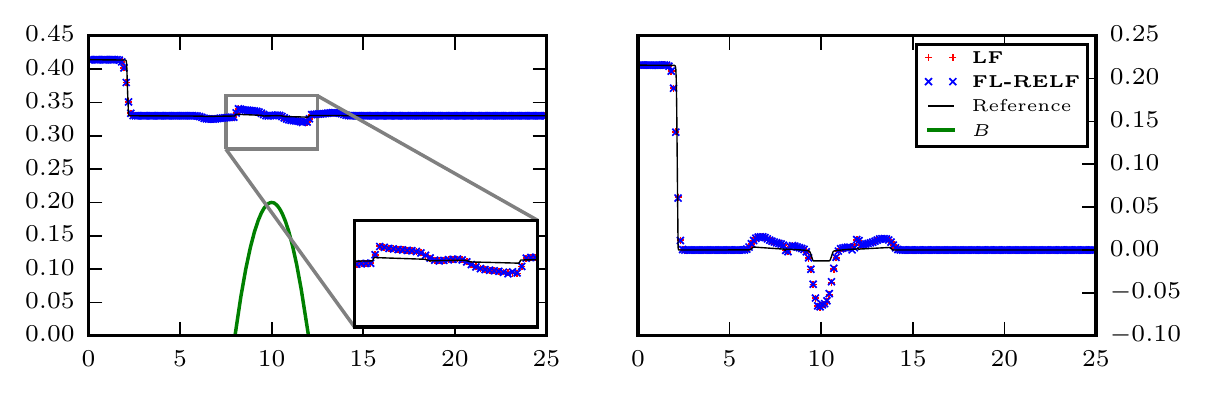} \\
            \includegraphics[scale=1]{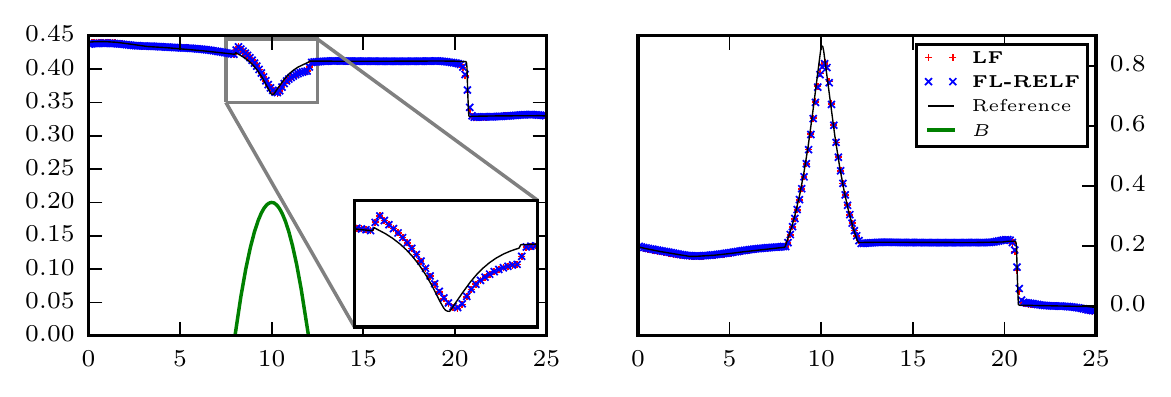} \\
            \includegraphics[scale=1]{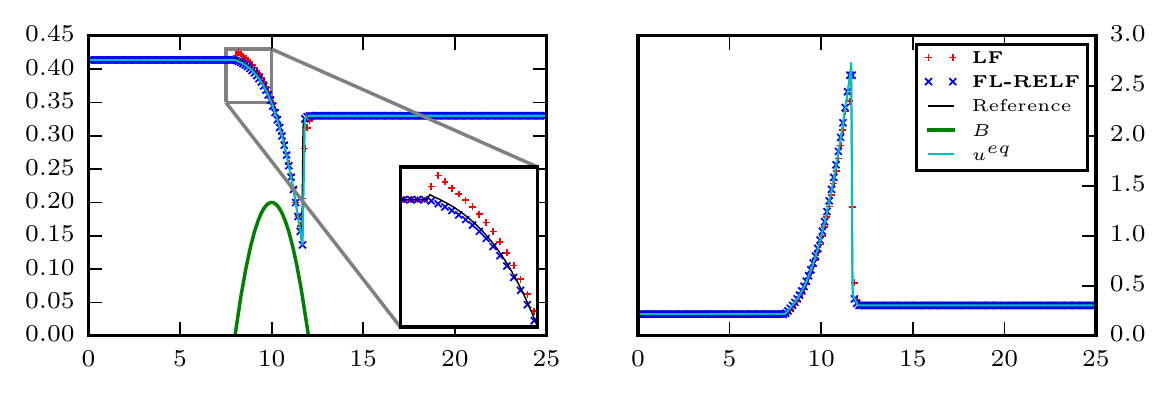}
          \end{tabular}
          \end{center}
          \caption{\textbf{Shallow water: transcritical flow.} Evolution of the water height $h$ (left) and Froude number $Fr= |v|/\sqrt{g h}$ (right) at times $t=1$, $10$ and $500$, for the classical (red dots) and flux limited equilibrium preserving (blue dashes) Lax-Friedrichs schemes. The green line represents the topography, the black dashes the reference solution and the cyan line, the equilibrium.}
          \label{fig:SW_transcriticalFL_Evol}
        \end{figure}  
        
        Nevertheless, one can overcome this by using the modified \emph{local residual equilibrium} scheme \eqref{eq:fl-re} (denoted by \textbf{FL-REFL}).
        Indeed, this local version can automatically switch between the classical and residual equilibrium schemes (we used in all the simulations the equilibrium indicator \eqref{def:eqIndic} given in the appendix), avoiding the spurious oscillations induced when the solution is very far from the equilibrium.
        Figure \ref{fig:SW_transcriticalFL_Evol} show the evolution of the solutions obtained with both \textbf{LF} and \textbf{FL-RELF} schemes, at times $t=1$, $10$ and $500$ (equilibrium), for $N_x = 200$ points in space. The reference solution has been obtained by the \textbf{LF} solver using $N_x = 1600$ points.
        We notice that the spurious spike has been damped by the local solver, remaining very close to the reference solution, and still resolve with an excellent accuracy the large time behavior of the solution. In particular, it avoids the wrong equilibrium solution given by the \textbf{LF} scheme.

        \begin{figure}
          \begin{center}
          \includegraphics[scale=.8]{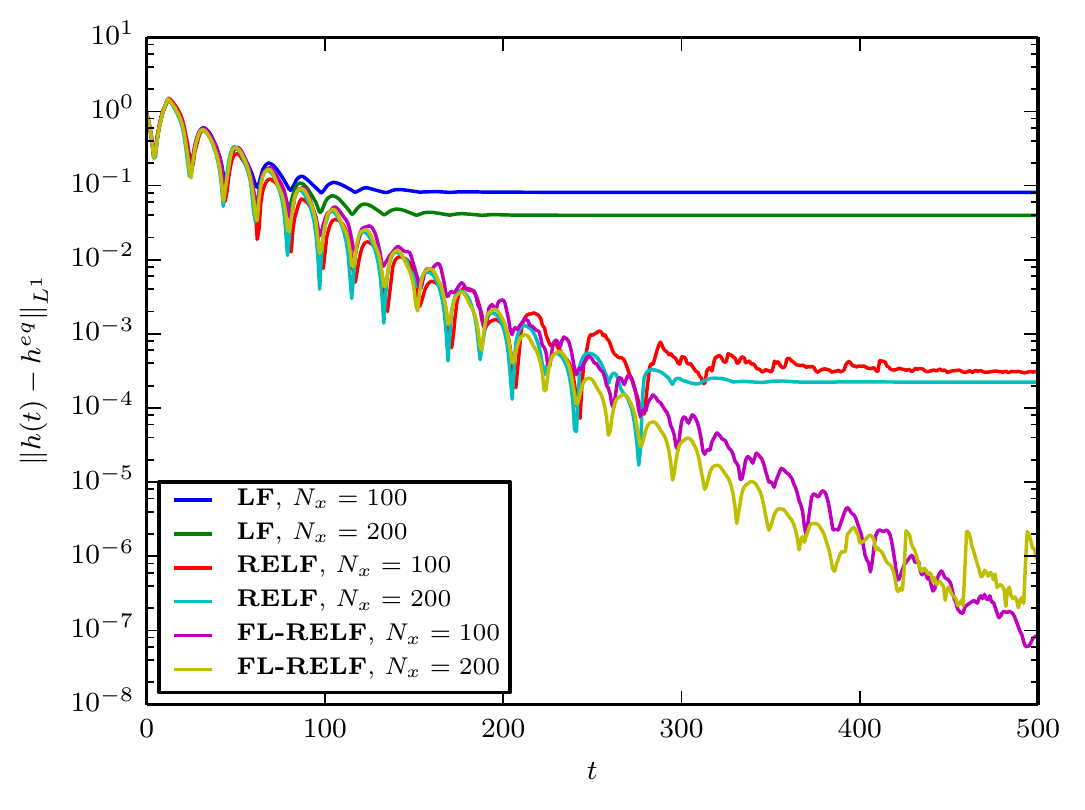}
          \caption{\textbf{Shallow water: transcritical flow.} Time evolution of the $L^1$ error of the surface water $h$, with respect to the equilibrium profile, for the standard (\textbf{LF}), residual equilibrium (\textbf{RELF}) and the local residual equilibrium (\textbf{FL-RELF}) Lax-Friedrichs schemes, with different mesh sizes.}
          \end{center}
          \label{fig:SW_Transcritical_error}          
        \end{figure}

        Finally, Figure \ref{fig:SW_Transcritical_error} presents the trend to equilibrium of the surface water $h$ towards its equilibrium value $h^{eq}$, for the three methods considered in this Section, namely the standard (\textbf{LF}), residual equilibrium (\textbf{RELF}) and the local residual equilibrium (\textbf{FL-RELF}) Lax-Friedrichs schemes.
        We observe that the residual equilibrium schemes greatly improve the large time behavior of the \textbf{LF} method, giving as expected a more accurate equilibrium profile. 
        Nevertheless, as we already noticed in Figures \ref{fig:SW_transcritical_Evol} and \ref{fig:SW_transcriticalFL_Evol}, the \textbf{RELF} needed some correction to resolve correctly the intermediate behavior of the solution, in particular because of the magnitude of the correction term \eqref{eq:re}.
        This is confirmed by this convergence rate, where the modified \textbf{FL-RELF} scheme is the only one able to approximate with a good accuracy the equilibrium state.
        
        We also notice in this Figure, thanks to the new schemes, that the convergence toward the equilibrium solution is exponential, with periodic oscillations. This looks related to the behavior observed for the solution of more general fluid equations, such that the Euler or Navier-Stokes systems.
        In these equations, this behavior can be explained by computing the dispersion relations of the equation, allowing roughly to show that the nonlinear semigroup of solution $\mathcal N_t$ behaves in large time as
        \[
          \left \| \mathcal{N}_t \right \| \leq C \exp(i \lambda_1 t - \lambda_2 t^2),
        \]
        where $\lambda_1 \in \RR$ and $\lambda_2 > 0$. 
        We can now  only conjecture that such a computation can also be done for this particular model, because we weren't able to prove anything about it.

%    \subsection{Power-law Tails}
%      \label{subNumPL}
%      $d = 1$, $\Omega = [-L,L]$,
%      
%      \paragraph{\textbf{Granular media with thermal bath.}}
%        
%        Equation:
%        \[
%          \frac{\partial u}{\partial t}(t,x) \, = \, \diverg \left( -\lambda W*u(y,x) + \beta x u(t,x( + \nabla u(t,x)\right ).
%        \]
%        
%        Equilibrium solution: thick Gaussians, but non explicit. Problems with the drift?
%      
%      
%      \paragraph{\textbf{Economy (Pareto tails).}}
%      
%        Non explicit equilibrium?

          \section{Conclusion}
          In this paper we introduced a general approach to deal numerically with a wide class of PDEs which admit steady states. The method is based on an underlying numerical discretization and is capable to preserve the steady states by keeping the accuracy of the underlying method. Therefore, the approach can be used also in conjunction with spectral techniques for which the direct construction of a steady state preserving approach is often impossible. Thanks to its structural simplicity it can be applied to very different PDEs, examples ranging from linear and nonlinear diffusion equations to Boltzmann equations and hyperbolic balance laws have been presented and show the effectiveness of the residual equilibrium approach.
          In principle the methods admits several improvement when applied to a specific kind of PDE, for example in order to preserve some peculiar property of the differential model like nonnegativity, monotonicity etc. Along this direction here we  discussed the case of the TVD property for hyperbolic balance laws using suitable equilibrium limiters in the scheme.

		\section*{Acknowledgment}          
		Thomas Rey would like to thanks Marianne Bessemoulin-Chatard for the many fruitful discussions they had about the residual equilibrium method and her help with the numerical methods for the porous medium equations.
		TR was partially funded by Labex CEMPI (ANR-11-LABX-0007-01).
		
   \appendix
   
          \section{Equilibrium limiters and TVD properties}
For clarity of approach here we consider the linear advection with source 
\begin{equation}
\frac{\partial u}{\partial t}=-a\frac{\partial u}{\partial x} -u, 
\label{eq:testp}
\end{equation}
where we assume $a>0$.

A conservative semi-discrete scheme has the general form 
\be
\frac{\partial u_i}{\partial t} =- a\frac{U_{i+1/2}-U_{i-1/2}}{\Delta x}-u_i,
\ee
where $U_{i\pm 1/2}$ are the edge fluxes for the $i$-th cell and 
\[
u_i=\frac1{\Delta x} \int_{x_{i-1/2}}^{x_{i+1/2}} u(x,t)\,dx,
\]
are the cell averages.

We denote by $u^{eq}$ the equilibrium state 
such that
\be
a\frac{\partial}{\partial x}u^{eq}+u^{eq}=0.
\label{eq:wbrel}
\ee
For example, if problem (\ref{eq:testp}) is considered for $x\in [0,+\infty)$ and supplemented with the boundary condition $u(0,t)=u_B$, $t>0$, 
we have explicitly
\[
u^{eq}(x)=u_B \exp\left\{-\frac{x}{a}\right\}.
\]
Let us denote by $U^{um}_{i\pm 1/2}$ the {numerical fluxes} of the underlying method. 
We define by $U^{eq}_{i\pm 1/2}$ the fluxes $U^{um}_{i\pm 1/2}$ evaluated by replacing $u$ with $u^{eq}$. Typically for these fluxes 
\[ 
a\frac{U^{eq}_{i+1/2}-U^{eq}_{i-1/2}}{\Delta x} - u^{eq}_i \neq 0.
\]
Now, we can define the {equilibrium preserving fluxes} 
\[
{\cal U}_{i\pm 1/2} = U^{um}_{i\pm 1/2}-U^{eq}_{i\pm 1/2},
\]
in order to construct the residual equilibrium semi-discrete scheme 
\be
\frac{\partial u_i}{\partial t} =- \frac{{\cal U}_{i+1/2}-{\cal U}_{i-1/2}}{\Delta x}+(u_i-u^{eq}_i).
\ee
To recover suitable monotonicity properties for the equilibrium preserving fluxes we add only a limited amount of the flux at equilibrium
\be
{\cal U}_{i\pm 1/2} = U^{um}_{i\pm 1/2}-\phi_{i} U^{eq}_{i\pm 1/2},
\ee
where $\phi_{i}=\phi(r_{i})\geq 0$ is an \emph{equilibrium flux limiter} and $r_{i}$ a suitable equilibrium indicator which depends on the solution. 

We now seek a flux limiter function and an equilibrium indicator in such a way that the scheme preserves simultaneously the steady states and some non-oscillatory property. Note that, since we are dealing with a linear problem, if the numerical scheme is also linear clearly the properties of the underlying method are kept by the residual equilibrium approach. In fact, the simple change of variable $u_i^n \to (u_i^n-u_i^{eq})$ makes the two schemes equivalent.

Therefore, we focus on nonlinear schemes and consider the explicit method
\be
u_i^{n+1} = u_i^n - \nu \Delta U^n_{i}-\Delta t\, u_i^n,
\label{eq:scheme}
\ee
where $\nu=a\Delta t/\Delta x$, $\Delta U^n_{i} = U^n_{i+1/2}-U^n_{i-1/2}$ and $U^n_{i\pm 1/2}$ is a second order TVD flux \cite{leveque:2002}. The TVD property corresponds to
\be
TV(u^{n+1})=\sum_i |u_{i+1}^{n+1}-u_i^{n+1}| \leq TV(u^n)=\sum_i |u_{i+1}^n-u_i^{n}|.
\ee
If we denote by 
\[
C_{i+1/2}=\nu \frac{\Delta U^n_{i+1}}{\Delta u^n_{i+1/2}},
\]
for scheme (\ref{eq:scheme}) it is immediate to show that the TVD property holds true if 
\be
0\leq C_{i+1/2}\leq 1-\Delta t.
\label{eq:tvd}
\ee
In fact we have
\begin{eqnarray*}
u_{i+1}^{n+1}-u_i^{n+1}&=&(u_{i+1}^{n}-u_i^{n})(1-\Delta t)-\nu (\Delta U^n_{i+1}-\Delta U^n_{i})\\
&=& (1-\Delta t)\Delta u^n_{i+1/2}-\nu\frac{\Delta U^n_{i+1}}{\Delta u^n_{i+1/2}}\Delta u^n_{i+1/2}+\nu\frac{\Delta U^n_{i}}{\Delta u^n_{i-1/2}}\Delta u^n_{i-1/2}.
\end{eqnarray*}
Under assumptions (\ref{eq:tvd}), summing up over all indexes we get
\begin{eqnarray*}
TV(u^{n+1}) &\leq & \sum_i (1-\Delta t-C_{i+1/2})|\Delta u^n_{i+1/2}|+\sum_i C_{i-1/2}|\Delta u^n_{i-1/2}|\\
&=& (1-\Delta t) \sum_i |\Delta u^n_{i+1/2}| =(1-\Delta t)\,TV(u^n).
\end{eqnarray*}
Finally, inequality (\ref{eq:tvd}) is recovered using suitable flux limiters \cite{Harten83, Sweby83}. 

Let us now include the equilibrium fluxes to obtain the scheme 
\bea
u_i^{n+1} &=& u_i^n - \nu\Delta U^n_{i} +\phi_{i}\nu\Delta U^{eq}_{i}-\Delta t\,\left(u_i^{n}- u_i^{eq}\right)\\
&=& u_i^n(1-\Delta t) - \nu\left(1-\phi_i\frac{\Delta U^{eq}_{i}}{\Delta U^n_{i}}\right)\Delta U^n_{i}+\Delta t\, u_i^{eq}.
\eea
Introducing the equilibrium indicator
\be
r_i=\frac{\Delta U^n_{i}}{\Delta U^{eq}_{i}},
\ee
we have
\[
u_i^{n+1}=u_i^n(1-\Delta t) - \nu\left(1-\frac{\phi(r_i)}{r_i}\right)\Delta U^n_{i}+{\Delta t}\, u_i^{eq}.
\]
To recover the TVD property we require that the equilibrium limiter function satisfies 
\be
0\leq \frac{\phi(r)}{r}\leq 1.
\label{eq:phib}
\ee
In fact, under this assumption, by the same arguments used for the underlying method we get
\[
TV(u^{n+1}) \leq (1-\Delta t)\, TV(u^n)+\Delta t\, TV(u^{eq}) \leq TV(u^n),
\]
provided that $TV(u^{eq}) \leq TV(u^n)$, $\forall\, n$.

If, in addition to (\ref{eq:phib}), we insist on $\phi(r)=0$ for $r \leq 0$ we have
\[
0\leq \phi(r) \leq r.
\]
Note, however, that we additionally want $\phi(r)\approx 1$ when $r$ is close to $1$ and $\phi(r)\approx 0$ when $r \gg 1$. For example, a possible choice is given by 
\begin{equation}
  \label{def:eqIndic}
\phi(r) = \begin{cases} r^{\alpha} &\mbox{if } 0< r \leq 1 \\ 
r^{-\alpha} & \mbox{if } r > 1, \end{cases} 
\end{equation}
where $\alpha > 1$ is a suitable parameter.

By  similar arguments as in \cite{Harten83, Sweby83} it is possible to extend the previous analysis to more general nonlinear problems of the type (\ref{eq:hypcon}) where the underlying method satisfies the TVD property and the source term $R(u)$ is such that $R'(u) \leq 0$.
\newpage

  \bibliographystyle{acm}
  \bibliography{biblioPR}

\end{document}